\documentclass[a4paper,12pt]{article}

\usepackage{latexsym}
\usepackage{mathrsfs}
\usepackage{amsfonts,amsmath,amssymb}
\usepackage{indentfirst}
\usepackage{subeqnarray}
\usepackage{verbatim}
\usepackage[pdftex]{graphicx}

\usepackage{epstopdf}
\usepackage{booktabs}
\usepackage{geometry}
\usepackage{bm}
\usepackage{enumerate}
\usepackage{float}
\usepackage{cite}
\usepackage{caption}
\usepackage{xcolor}
\usepackage{enumitem}
\makeatletter
\newcommand{\itemlabel}[2]{#2\def\@currentlabel{#2}\label{#1}}
\makeatother

\usepackage{algorithmicx}
\usepackage[ruled]{algorithm}
\usepackage{algpseudocode}

\def\real{\mathbb{R}}

\newtheorem{theorem}{Theorem}
\newtheorem{lemma}[theorem]{Lemma}

\numberwithin{equation}{section}
\newenvironment{pfs}[1][Proof]{\noindent\textbf{#1.} }
{\ \rule{0.75em}{0.75em}\smallskip}

\DeclareMathOperator*{\argmin}{arg\,min}

\begin{document}


\begin{center}
 \Large\bf A survey of numerical methods \\ for hemivariational inequalities \\ with applications to Contact Mechanics
\end{center}

\begin{center}
Anna Ochal\footnote{Jagiellonian University in Krakow, Faculty of Mathematics and Computer Science, 
Lojasiewicza 6, 30-348 Krakow, Poland. Email: {\tt anna.ochal@uj.edu.pl}}, 
\quad \quad
Michal Jureczka\footnote{Jagiellonian University in Krakow, Faculty of Mathematics and Computer Science, 
Lojasiewicza 6, 30-348 Krakow, Poland. Email: {\tt michal.jureczka@uj.edu.pl}}\quad and \quad 
Piotr Bartman\footnote{Jagiellonian University in Krakow, Faculty of Mathematics and Computer Science, 
Lojasiewicza 6, 30-348 Krakow, Poland. Email: {\tt piotr.bartman@doctoral.uj.edu.pl}}
\end{center}

\begin{quote}
{\bf Abstract}. 
In this paper we present an abstract nonsmooth optimization problem for which we recall existence and uniqueness results. 
We show a numerical scheme to approximate its solution. 
The theory is later applied to a sample static contact problem describing an elastic body in frictional contact with a foundation. 
This problem leads to a hemivariational inequality which we solve numerically.
Finally, we compare three computational methods of solving contact mechanical problems: direct optimization method, augmented Lagrangian method and primal-dual active set strategy.



\end{quote}

\begin{center}
{\it Dedicated to 60-th birthday of Professor Stanisław Migórski}
\end{center}


\section{Introduction}
\label{sec:intro}

\noindent
Mathematical models which describe contact between a deformable body and a foundation have various applications. 
Many of them have already been analyzed in the literature, where behavior of the body on the contact boundary is governed by monotone functions responsible in turn for foundation response in the normal direction to the contact boundary and friction in the tangent plane to the boundary. 
However, considering nonmonotone functions requires a different analytical as well as numerical treatment.
For example, let us consider a foundation made from several layers with different properties, so we have to consider different friction laws along the penetration. 
This leads to contact mechanical problem which involves nonmonotone functions.

To find an approximate solution, we first formulate an abstract scheme for the chosen class of contact mechanical problem. 
It starts with an introduction of a general nonsmooth optimization problem together with required assumptions followed by existence and uniqueness results. 
Numerical approximation of optimization problem with obtained error estimate allows us to use this abstract scheme to a static contact problem.
In this paper we consider a nonmonotone friction law that depends on both normal and tangential components of displacement and its weak formulation which leads to hemivariational inequality. 
We consider a similar mechanical model to the one described in \cite{JO}, but with a more general $A$ operator as potential operator.
Next, we compare three popular methods of solving the introduced problem: direct optimization method, augmented Lagrangian method and primal-dual active set strategy.

The definition and properties of the Clarke subdifferential and tools used to solve optimization problems can be found in \cite{C}, differences between nonsmooth and nonconvex optimization methods in \cite{BKM}, and details on computational contact mechanics in \cite{W}. 
Introduction to the theory of hemivariational inequalities is available in \cite{P}, and the first usage of the finite element method to solve these inequalities is in \cite{HMP}.
Early study of vector-valued hemivariational problems related to FEM is presented in \cite{MH} and recent analysis of hemivariational and variational-hemivariational inequalities was presented in \cite{MOS2, MOS}. Numerical analysis of such problems can be found for example in papers \cite{BBK, BBKR, BHM, H, HSB, HSD}.

In \cite{H} is presented an error estimate of stationary variational-hemivariational inequalities. 
In our paper variational part of inequality is not present and the inequality is not constrained, nevertheless to reflect the dependence of friction law on the normal component of the displacement error estimate had to be generalized.

The direct optimization method was previously compared to the augmented Lagrangian method in \cite{BBK}. Nevertheless, usually in papers containing Contact Mechanics simulations one method is chosen and presented. Early ideas about optimization in Contact Mechanics were presented in \cite{P} and further details can be found in \cite{JO, P}. 
The augmented Lagrangian method was reviewed or applied in \cite{AC, BBKR, PC, W}. In this paper we also include a third method, called primal-dual active set strategy, for comparison of implementation and obtained results. Details on this third method can be found in \cite{K, BCKYZ, XCHQ}.

The paper is organized as follows. 
In Section~\ref{sec:op}, we formulate a general differential inclusion and equivalent optimization problem. 
We present existence and uniqueness results under usual assumptions.
Next, we introduce the discrete formulation of the optimization problem together with numerical error estimate theorem.
The introduced abstract scheme is then used for a weak formulation of chosen contact mechanical problem and presented in Section~~\ref{sec:mcp}.
In Section~\ref{sec:mo}, we briefly describe three alternative methods of solving problem from the previous section: direct optimization method, augmented Lagrangian method and primal-dual active set strategy. 
Finally, we compare the results of the error estimate obtained for each method.

\section{A general optimization problem} \label{sec:op}
\noindent
In this section we recall notation, definitions and preliminary material (we refer \cite{C, HMP, Z}), and to analysis of a general optimization problem.
For a normed space $X$, we denote by $\|\cdot\|_X$ its norm, by $X^*$ its topological dual and by $\langle\cdot,\cdot\rangle_{X^*\times X}$ the duality pairing of $X^*$ and $X$. 
Given two normed spaces X and Y, $\mathcal{L}(X, Y)$ is the space of all linear continuous operators from $X$ to $Y$ with the norm $\|\cdot\|_{\mathcal{L}(X, Y)}$.
Let $L \in \mathcal{L}(X, Y)$, then the adjoint operator to $L$ is denoted by $L^* \colon Y^* \to X^*$.
Let $X$ be a real Banach space, and let $j \colon X \to \real$ be locally Lipschitz continuous. 
Then the generalized (Clarke) directional derivative of $j$ at $x \in X$ in the direction $v \in X$ is
\begin{align*}
  &j^0(x;v) := \limsup_{y \to x, \lambda \searrow 0} \frac{j(y+\lambda v)  - j(y)}{\lambda}.
\end{align*}
The generalized subdifferential of $j$ at $x$ is
\begin{align*}
  &\partial j(x) := \{\xi \in X^*\, | \, \langle \xi, v\rangle_{X^*\times X}\leq j^0(x;v) \  \mbox{ for all } v \in X \}.
\end{align*}
If $\partial j(x)$ is nonempty, then any element $\xi \in \partial j(x)$ is called a subgradient of $j$ at $x$ (cf. \cite{C}).
If $j\colon X^n \to \mathbb{R}$ is a locally Lipschitz function of $n$ variables, then we use $\partial_i j$ and $j_i^0$ to denote the Clarke subdifferential and generalized directional derivative with respect to $i$-th variable of~$j$, respectively.

Recall that an operator $A \colon X \to X^* $ is called a potential operator if there exists a~G\^ateaux differentiable functional $ F_A \colon X \to \mathbb{R}$ such that $ A = F_A^{\prime} $.
The functional $F_A$ is called a potential of $A$.
Basic properties of potential operators can be found in \cite{Z}.
We recall that $A \in \mathcal{L}(X, X^*)$ is a potential operator if and only if it is symmetric.
Moreover, under this symmetry condition, a potential functional is given by
\begin{align*}
  &F_A(v) = \frac{1}{2} \langle Av, v\rangle_{X^*\times X} \  \mbox{ for all } v \in X.
\end{align*}
Throughout the paper, by $c>0$ we denote a generic constant whose value may change from one place to another but it is independent of other quantities of concern.

Let now $V$ be a reflexive Banach space and $X$ be a Banach space. 
Given an operator $A \colon V \to V^*$, a locally Lipschitz function $J \colon X \times X \to \real$, a linear operator $\gamma \colon V \to X$, and a~linear functional $f \colon V \to \real$, we consider the following operator inclusion problem

\medskip
\noindent
\textbf{Problem $\bm{P_{incl}}$:} {\it\ Find $u \in V$ such that}
\begin{align*}
  &A u + \gamma^* \partial_2 J(\gamma u, \gamma u) \ni f. 
\end{align*}

\noindent
We say that $u\in V$ is a solution to Problem~$P_{incl}$ if there exists $z\in \partial_2 J(\gamma u, \gamma u)$ such that $A u + \gamma^* z =f$. 

In the study of Problem $P_{incl}$ we adopt the following hypotheses

\begin{enumerate}[wide, labelwidth=0pt, labelindent=0pt]
    \item[\itemlabel{assumption:A}{$H(A)$}]:
    \quad The operator $A \colon V \to V^*$ is such that
    \begin{enumerate}[label=(\alph*)]
      \item \label{assumption:A:Lipschitz_continuous}
         $A$ is Lipschitz continuous, i.e.,  $\|A u - A v \|_{V^*} \leq L_A \, \|u - v \|_{V}$ for all $u, v \in V$ with $L_A > 0$,
      \item \label{assumption:A:potential_operator}
         $A$ is a potential operator with a potential~$F_A$,
      \item \label{assumption:A:strongly_monotone}
        $A$ is strongly monotone, i.e., 
         $\langle A u - A v, u - v\rangle_{V^*\times V} \geq m_A \|u - v\|_V^2$ for all $u, v \in V$ with $m_A>0$.
    \end{enumerate}

    \item[\itemlabel{assumption:J}{$H(J)$}]:
    \quad The functional $J \colon X \times X \to \mathbb{R}$ satisfies
    \noindent
    \begin{enumerate}[label=(\alph*)]
      \item \label{assumption:J:lipschitz}
         $J$ is locally Lipschitz continuous with respect to its second variable,
      \item \label{assumption:J:second}
         there exist $c_{0}, c_{1}, c_{2} \geq 0$ such that \\[2mm]
         \hspace*{1cm}$\|\partial_2 J(w, v)\|_{X^*} \leq c_{0} + c_{1}\| v \|_X + c_{2}\| w \|_X$\quad for all $w, v \in X$,
      \item \label{assumption:J:third}
        there exist $m_\alpha, m_L  \geq 0$ such that \\[2mm]
        $J_2^0(w_1, v_1; v_2 - v_1) + J_2^0(w _2, v_2; v_1 - v_2)\leq m_\alpha\| v_1- v_2\|_X^2+ m_L\| w_1 -  w_2 \|_{X} \| v_1 -  v_2\|_{X}$ \\[2mm]
        for all $ w_1,  w_2,  v_1,  v_2 \in X$.
    \end{enumerate}

    \item[\itemlabel{assumption:f}{$H(\gamma, f)$}]:
    \quad $\gamma \in \mathcal{L}(V, X)$, \quad $f \in V^*$.
    
    \item[\itemlabel{assumption:smallness}{$(H_s)$}]:
    \quad $m_A > (m_\alpha + m_L)\, c_\gamma^2$, where $c_\gamma := \|\gamma \|_{\mathcal{L}(V, X)}$.
\end{enumerate} 

\bigskip
It is easy to see that \ref{assumption:J}\ref{assumption:J:third} is equivalent to the following condition
\begin{align*}
 &\langle \partial_2 J(w_1,  v_1) - \partial_2 J(w_2,  v_2), v_1 -  v_2 \rangle_{X^*\times X} \geq - m_\alpha\| v_1- v_2\|_X^2 - m_L\| w_1- w_2\|_X \| v_1- v_2\|_X
\end{align*}
for all  $w_1, w_2, v_1, v_2 \in X$. 
We remark that this condition generates the relaxed monotonicity condition which holds in a case of $J$ independent of its first variable, i.e.,
\begin{align}
 &\langle \partial J(v_1) - \partial J(v_2), v_1 -  v_2 \rangle_{X^*\times X} \geq - m_\alpha \| v_1- v_2\|_X^2.  \label{eq:RM}
\end{align}
for all  $v_1, v_2 \in X$.

Under introduced assumptions, we consider the following optimization problem

\medskip
\noindent
\textbf{Problem $\bm{P_{opt}}$:} {\it \ Find $ u  \in V$ such that}
\begin{align*}
  0 \in \partial_2 \mathcal{L}( u , u ).
\end{align*}
Here,  $\mathcal{L}: V \times V \rightarrow \mathbb{R}$ is defined by
\begin{align}
  \mathcal{L}( w ,  v ) = F_A(v)  - \langle  f ,  v  \rangle_{V^*\times V} + J(\gamma  w  ,\gamma  v  ). \label{eq:L}
\end{align}
for all $ w ,  v  \in V$.

We start with recalling some properties of the functional~$\mathcal{L}$.
This result is followed by argument similar to the ones used in \cite[Lemma 2]{JO}  but with more general operator $A$ being a potential operator (see also \cite[Proposition 2.5]{Han}). 
\begin{lemma} \label{lem:Lprop}
If the hypotheses \ref{assumption:A}, \ref{assumption:J}, \ref{assumption:f} and \ref{assumption:smallness} hold, then for a fixed $w\in V$ the functional  $\mathcal{L}(w,\cdot)\colon V \rightarrow \mathbb{R}$, defined by {\rm{(\ref{eq:L})}},
is locally Lipschitz continuous and strictly convex, hence also coercive, and 
$$\partial_2 \mathcal{L}( w ,  v ) \subseteq A  v  -  f  + \gamma^* \partial_2 J(\gamma  w , \gamma  v ).$$
\end{lemma}

\medskip
\noindent

The following result shows the relation between Problems $P_{incl}$ and $P_{opt}$ as well the existence of their unique solution.

\begin{theorem} \label{I=O}
If the hypotheses $H(A)$, $H(J)$, \ref{assumption:f} and $(H_s)$ hold, then Problems~$P_{incl}$ and $P_{opt}$ are equivalent,
they have a unique solution $u\in V$ and this solution satisfies estimation
\begin{align} \label{eq:No2}
  &\| u \|_V \leq c(1+\| f \|_{V^*})
\end{align}
with a positive constant $c$.
\end{theorem}

\noindent
Detailed arguments can be found in~\cite{JO} and are omitted here. We only mention the main steps of the proof. We first observe that Lemma~$\ref{lem:Lprop}$ implies that every solution to Problem $P_{opt}$ solves Problem~$P_{incl}$. 
Moreover, if Problem~$P_{incl}$ has a solution, then it is unique.
Then, it can be shown that Problem~$P_{opt}$ has a unique solution. This follows from the Banach fixed point theorem applied to an operator $\Lambda \colon V \to V$ given by
\begin{align*}
  & \Lambda w := \argmin_{ v \in V} \mathcal{L}(w, v)  \  \mbox{ for all } w \in V.
\end{align*}
Using the above facts, we see that a unique solution to Problem~$P_{opt}$ is also a unique solution to Problem~$P_{incl}$. 
And because of the uniqueness of the solution to Problem~$P_{incl}$ we deduce that Problem~$P_{incl}$ and Problem~$P_{opt}$ are equivalent. Finally, the inequality (\ref{eq:No2}) is a consequence of 
\ref{assumption:J}\ref{assumption:J:second}-\ref{assumption:J:third}, \ref{assumption:A}\ref{assumption:A:strongly_monotone} and \ref{assumption:smallness}. 

\medskip

We are now in a position to present numerical methods for solving the optimization problem. We keep assumptions $H(A)$, $H(J)$, $H(\gamma, f)$ and $(H_s)$ so that Problem~$P_{opt}$ has a unique solution~$u\in V$.
Let $V^h \subset V$ be a finite dimensional subspace with a discretization parameter $h>0$. We consider the following discrete scheme of Problem~$P_{opt}$.

\bigskip

\noindent
\textbf{Problem $\bm{P_{opt}^{h}}$:} {\it \ Find  $u^h \in V^h$ such that}
\begin{align*}
0 \in \partial_2 \mathcal{L}(u^h,  u^h).
\end{align*}
We can apply the arguments of the proof of Theorem~$\ref{I=O}$ in the setting of the finite dimensional space~$V^h$, to conclude the existence of a unique solution to Problem~$P_{opt}^h$ and equivalence to the discrete version of Problem~$P_{incl}$.
We now present the theorem concerning the error estimate of the introduced numerical scheme.

\begin{theorem} \label{estimate}
If the hypotheses $H(A)$, $H(J)$, $H(\gamma, f)$ and $(H_s)$ hold, then for the unique solutions $u$ and $u^h$ to Problems~$P_{opt}$ and $P_{opt}^h$, respectively,  there exists a constant $c>0$ such that
\begin{equation} \label{thminequality}
  \|  u  -  u ^h\|_V^2 \leq c\,\inf\limits_{ v ^h \in V^h}  \Big\{ \| u  -  v ^h \|_V^2 + \|\gamma  u  - \gamma v ^h \|_X + R( u ,  v ^h) \Big\},
\end{equation}
 where a residual quantity is defined by
\begin{equation} \label{R}
  R( u ,  v ^h) =  \langle A  u ,  v ^h -  u  \rangle_{V^*\times V} + \langle  f ,  u  -  v ^h \rangle_{V^*\times V}.
\end{equation}
\end{theorem}

\begin{pfs}
\noindent
Let $u$ and $u^h$ be solutions to Problems~$P_{opt}$ and $P_{opt}^h$, respectively. 
Hence, they satisfy the corresponding inclusion problems and the following inequalities, respectively
\begin{align*}
&\langle f - Au, v  \rangle_{V^*\times V}  \leq J_2^0(\gamma u, \gamma u; \gamma v) \quad \mbox{for \ all\ } v \in V, \\[2mm] 
&\langle f - Au^h, v  \rangle_{V^*\times V}  \leq J_2^0(\gamma u^h, \gamma u^h; \gamma v) \quad \mbox{for \ all\ } v \in V^h. 
\end{align*}
Setting $v=u^h-u$ in the first inequality, and $v=v^h-u^h$ with $v^h \in V^h$ in the second one, then adding the resulting inequalities, we deduce for all $ v^h \in V^h$
\begin{align*}
  &\langle f, v ^h - u \rangle_{V^*\times V} + \langle Au^h - Au, u^h - u \rangle_{V^*\times V} - \langle A u ^h, v ^h -  u  \rangle_{V^*\times V} \nonumber\\
  &\qquad \leq J_2^0(\gamma u , \gamma u ; \gamma u ^h - \gamma u ) + J_2^0(\gamma  u ^h, \gamma  u ^h; \gamma  v ^h - \gamma  u ^h). \label{27}
\end{align*}
The subadditivity of generalized directional derivative (cf. \cite{MOS}) and $H(J)$(c), give
\begin{align*}
  & J_2^0(\gamma u, \gamma u; \gamma u ^h - \gamma u) + J_2^0(\gamma u^h, \gamma u^h; \gamma v^h - \gamma u^h )\nonumber\\
  &\leq J_2^0(\gamma u, \gamma u; \gamma u^h - \gamma u) + J_2^0(\gamma u^h, \gamma u^h; \gamma u - \gamma u^h) + J_2^0(\gamma u^h, \gamma u^h; \gamma v^h - \gamma u)\nonumber\\
  &\leq (m_\alpha + m_L) \|\gamma u^h -\gamma u \|_X^2 + \left( c_0 + (c_1+c_2) \|\gamma u^h\|_X \right)\| \gamma  v ^h - \gamma  u \|_X. 
\end{align*}
From Theorem~$\ref{I=O}$ applied to discrete version of Problem $P_{incl}$ we obtain the uniform boundedness property with respect to $h$
\begin{equation*}
\|\gamma u^h\|_X \leq c_{\gamma}\|u^h\|_V\leq c\,(1+\|f\|_{V^*}).
\end{equation*}
Hence, combining the above inequalities, we deduce for all $ v ^h \in V^h$
 \begin{align*}
  &\langle Au^h - Au, u^h - u \rangle_{V^*\times V} \leq  \langle A u ^h - Au, v ^h -  u  \rangle_{V^*\times V} + \langle A u, v ^h -  u  \rangle_{V^*\times V} \nonumber\\
  &\qquad + \langle f, u-v^h \rangle_{V^*\times V} + (m_\alpha + m_L)c^2_\gamma \|u^h - u \|_V^2 + c \, \|\gamma v^h  - \gamma u\|_X .
\end{align*}
Using definition~(\ref{R}) and assumption~$H(A)$, we have for all $v^h\in V^h$
\begin{align*}
  &m_A \|  u ^h -  u \|_V^2 \leq L_A \, \| u ^h -  u \|_V\| v ^h -  u \|_V + R( u ,  v ^h) \nonumber\\
   &\qquad +  (m_\alpha + m_L) c_\gamma^2\|  u  -  u ^h \|_V^2 + c \, \| \gamma  u  - \gamma  v ^h \|_X. 
\end{align*}
Finally, the Cauchy inequality with $\varepsilon > 0$ yields
\begin{align*}
  & m_A \| u  -  u ^h\|_V^2 \leq \varepsilon\|  u  -  u ^h\|_V^2 + \frac{L_A^2}{4\varepsilon}\|  u  -  v ^h\|_V^2 + R( u ,  v ^h) \\
  &\qquad + (m_\alpha + m_L) c_\gamma^2\|  u  -  u ^h \|_V^2 + c \, \| \gamma  u  - \gamma  v ^h \|_X
\end{align*}
which implies for all $ v ^h \in V^h$
\begin{align*}
  & \Big(m_A - (m_\alpha + m_L) c_\gamma^2 - \varepsilon\Big) \|  u  -  u ^h\|_V^2 \leq \frac{c}{\varepsilon}\| u  -  v ^h \|_V^2 + R( u ,  v ^h) + c \, \| \gamma  u  - \gamma  v ^h \|_X.
\end{align*}
For sufficiently small $\varepsilon$ and by $(H_s)$, we obtain the desired C\'ea type inequality.
\end{pfs}

\section{Application to Contact Mechanics} \label{sec:mcp}
\noindent

This section presents a sample mechanical contact problem where results of the previous section are applied.
We want to find the body displacement in a static state.
At the beginning we introduce the physical setting and notation.

Let us consider an elastic body in a domain $\Omega \subset \mathbb{R}^{d}$, where $d = 2, 3$ in application. 
Boundary of $\Omega$ is denoted as $\Gamma$ and is divided into three disjoint measurable parts:
$\Gamma_{C}, \Gamma_{N}, \Gamma_{D}$, where the measure of part $\Gamma_D$ is positive.
Moreover $\Gamma$ is Lipschitz continuous, so the outward normal vector $\bm{\nu}$ to $\Gamma$ exists a.e. on the boundary.
To model contact with the foundation on boundary $\Gamma_{C}$ we use general subdifferential inclusions. 
Displacement of the body is equal $\bm{0}$ on $\Gamma_{D}$,
a surface force of density $\bm{f}_N$ acts on the boundary~$\Gamma_{N}$ and a body force of density $\bm{f}_0$ acts in $\Omega$. 

Let us denote by ``$\cdot$'' and $\|\cdot\|$ the scalar product and the Euclidean norm in $\mathbb{R}^{d}$ or $\mathbb{S}^{d}$, respectively, where  $\mathbb{S}^{d} = \mathbb{R}^{d \times d}_{sym}$. Indices $i$ and $j$ run from $1$ to $d$ and summation over repeated indices is implied. We denote the divergence operator by $\textrm{Div }\bm{\sigma} = \left(\frac{\partial\sigma_{ij}}{\partial x_j}\right)$. 
The linearized (small) strain tensor for displacement $\bm{u} \in H^1(\Omega)^d$ is defined by
\begin{equation}\nonumber
 \bm{\varepsilon}(\bm{u})=(\varepsilon_{ij}(\bm{u})), \quad \varepsilon_{ij}(\bm{u}) = \frac{1}{2}\left(\frac{\partial u_{i}}{\partial x_j} + \frac{\partial u_{j}}{\partial x_i}\right).
\end{equation}
Let $u_\nu= \bm{u}\cdot \bm{\nu}$ and $\sigma_\nu= \bm{\sigma}\bm{\nu} \cdot \bm{\nu}$ be the normal components of $\bm{u}$ and $\bm{\sigma}$, respectively,  and let $\bm{u}_\tau =\bm{u}-u_\nu\bm{\nu}$ and $\bm{\sigma}_\tau =\bm{\sigma}\bm{\nu}-\sigma_\nu\bm{\nu}$ be their tangential components, respectively. In what follows, for simplicity, we sometimes do not indicate explicitly the dependence of various functions on the spatial variable $\bm{x}$.

Now let us introduce the classical formulation of the considered mechanical contact problem.

\medskip
\noindent
\textbf{Problem $\bm{P}$:} \textit{Find a displacement field $\bm{u}\colon \Omega \rightarrow \mathbb{R}^{d}$ and a stress field $\bm{\sigma}\colon \Omega \rightarrow \mathbb{S}^{d}$ such that}

\begin{align}
  \bm{\sigma}  = \mathcal{A}(\bm{\varepsilon}(\bm{u})) \qquad &\textrm{ in } \Omega \label{P1}\\
  \textrm{Div }\bm{\sigma} + \bm{f}_{0} = \bm{0}  \qquad &\textrm{ in } \Omega \label{P2}\\
  \bm{u} = \bm{0}  \qquad &\textrm{ on } \Gamma_{D} \label{P3}\\
  \bm{\sigma}\bm{\nu} = \bm{f}_{N} \qquad &\textrm{ on } \Gamma_{N} \label{P4}\\
   -\sigma_{\nu} \in \partial j_{\nu}(u_{\nu}) \qquad &\textrm{ on } \Gamma_{C} \label{P5}\\
   -\bm{\sigma_{\tau}} \in h_{\tau}\,\partial j_{\tau}(\bm{u_{\tau}}) \qquad &\textrm{ on } \Gamma_{C} \label{P6}
\end{align}

\noindent
Here, equation~\eqref{P1} represents an elastic constitutive law and $\mathcal{A}$ is an elasticity operator.
Equilibrium equation~(\ref{P2}) reflects the fact that the problem is static.
Equation~(\ref{P3}) represents the clamped boundary condition on $\Gamma_{D}$ and~(\ref{P4}) represents tractions applied on $\Gamma_{N}$.
Inclusion~(\ref{P5}) describes the response of the foundation in normal direction, whereas the friction is modeled by inclusion~(\ref{P6}), where $j_\nu$ and $j_\tau$ are given superpotentials, and $h_\tau$ is a~given friction bound. Note that to simplify simulation, the function $h_\tau$ does not depend on $u_\nu$.


\medskip
\noindent
Now we present the hypotheses on data of Problem~$P$.

\medskip
\noindent
\begin{enumerate}[wide, labelwidth=0pt, labelindent=0pt]
    \item[\itemlabel{assumption:$A$}{$H({\mathcal{A}})$}]:
    \quad ${\mathcal{A}} \colon \Omega \times {\mathbb S}^d \to {\mathbb S}^d$ satisfies
    \begin{enumerate}[label=(\alph*)]
      \item \label{assumption:$A$:a}
        $\mathcal{A}(\bm{x},\bm{\tau}) = (a_{ijkh}(\bm{x})\tau_{kh})$
        for all $\bm{\tau} \in {\mathbb S}^d$, a.e. $\bm{x}\in\Omega,\ a_{ijkh} \in L^{\infty}(\Omega),$
      \item \label{assumption:$A$:b}
        $\mathcal{A}(\bm{x},\bm{\tau}_1) \cdot \bm{\tau}_2 = \bm{\tau}_1 \cdot  \mathcal{A}(\bm{x},\bm{\tau}_2)$ for all $\bm{\tau}_1, \bm{\tau}_2 \in {\mathbb S}^d$, a.e. $\bm{x}\in\Omega$,      
      \item \label{assumption:$A$:c}
        there exists $m_{\mathcal{A}}>0$ such that $\mathcal{A}(\bm{x},\bm{\tau}) \cdot \bm{\tau} \geq m_{\mathcal{A}} \|\bm{\tau}\|^2$ for all $\bm{\tau} \in {\mathbb S}^d$, a.e. $\bm{x}\in\Omega$.
     \end{enumerate}
     
    \item[\itemlabel{assumption:j_nu}{$H(j_{\nu})$}]:
    \quad $j_{\nu} \colon \Gamma_C \times \mathbb{R} \to \mathbb{R}$ satisfies
    \begin{enumerate}[label=(\alph*)]
      \item \label{assumption:j_nu:a}
         $j_{\nu}(\cdot, \xi)$ is measurable on $\Gamma_C$ for all $\xi \in \mathbb{R}$ and there exists $e \in L^2(\Gamma_C)$ such that \\ $j_{\nu}(\cdot,e(\cdot))\in L^1(\Gamma_C)$,
      \item \label{assumption:j_nu:b}
        $j_{\nu}(\bm{x}, \cdot)$ is locally Lipschitz continuous on $\mathbb{R}$ for a.e. $\bm{x} \in \Gamma_C$,
      \item \label{assumption:j_nu:c}
        there exist $c_{\nu0}, c_{\nu1} \geq 0$ such that \\[2mm]
        \hspace*{1cm}$|\partial_2 j_{\nu}(\bm{x}, \xi)| \leq c_{\nu0} + c_{\nu1}|\xi|$\quad for all $\xi \in \mathbb{R}$, a.e. $\bm{x} \in \Gamma_C$,
      \item \label{assumption:j_nu:d}
        there exists $\alpha_{\nu} \geq 0$ such that \\[2mm]
        \hspace*{1cm}$(j_{\nu})_2^0(\bm{x},\xi_1;\xi_2-\xi_1) + (j_{\nu})_2^0(\bm{x},\xi_2;\xi_1-\xi_2)\leq \alpha_{\nu}|\xi_1-\xi_2|^2$ \\[2mm]
     for all $\xi_1, \xi_2 \in \mathbb{R}$, a.e. $\bm{x} \in \Gamma_C$.
     \end{enumerate}
     
    \item[\itemlabel{assumption:j_tau}{$H(j_{\tau})$}]:
    \quad $j_{\tau} \colon \Gamma_C \times \mathbb{R}^{d} \to \mathbb{R}$ satisfies
    \begin{enumerate}[label=(\alph*)]
      \item \label{assumption:j_tau:a}
        $j_{\tau}(\cdot, \bm{\xi})$ is measurable on $\Gamma_C$ for all $\bm{\xi} \in \mathbb{R}^{d}$ and there exists $\bm{e} \in L^2(\Gamma_C)^{d}$ such that $j_{\tau}(\cdot,\bm{e}(\cdot))\in L^1(\Gamma_C)$,
      \item \label{assumption:j_tau:b}
        there exists $c_{\tau}>0$ such that \\[2mm]
        \hspace*{1cm}$|j_{\tau}(\bm{x}, \bm{\xi}_1) - j_{\tau}(\bm{x}, \bm{\xi}_2)| \leq c_{\tau} \|\bm{\xi}_1 - \bm{\xi}_2\|$\quad for all $\bm{\xi}_1, \bm{\xi}_2 \in \mathbb{R}^d$, a.e. $\bm{x} \in \Gamma_C$,
      \item \label{assumption:j_tau:c}
        there exists $\alpha_{\tau} \geq 0$ such that  \\[2mm]
        \hspace*{1cm} $(j_{\tau})_2^0(\bm{x},\bm{\xi}_1;\bm{\xi}_2-\bm{\xi}_1) + (j_{\tau})_2^0(\bm{x},\bm{\xi}_2;\bm{\xi}_1-\bm{\xi}_2)\leq \alpha_\tau\|\bm{\xi}_1-\bm{\xi}_2\|^2$\\[2mm]
        for all $\bm{\xi}_1, \bm{\xi}_2 \in \mathbb{R}^{d}$, a.e. $\bm{x} \in \Gamma_C$.
     \end{enumerate}
     
    \item[\itemlabel{assumption:h}{$H(h)$}]:
    \quad $h_{\tau} \colon \Gamma_C \to \mathbb{R}$ satisfies
    \begin{enumerate}[label=(\alph*)]
      \item \label{assumption:h:a}
        $h_{\tau}(\cdot)$ is measurable on $\Gamma_C$
      \item \label{assumption:h:b}
        there exists $\overline{h}_{\tau} > 0$ such that $0 \leq h_{\tau}(\bm{x}) \leq \overline{h}_{\tau}$ a.e. $\bm{x} \in \Gamma_C$,
     \end{enumerate}
     
    \item[\itemlabel{assumption:H_0}{$(H_0)$}]:
      \quad $\bm{f}_0 \in L^2(\Omega)^d, \quad \bm{f}_N \in L^2(\Gamma_N)^d$.
\end{enumerate}

\bigskip
Note that condition $H(j_{\tau})$(b) is equivalent to the fact that $j_{\tau}(\bm{x},\cdot)$ is locally Lipschitz continuous and there exists $c_{\tau} \geq 0$ such that $\|\partial_2 j_{\tau}(\bm{x} ,  \bm{\xi} )\| \leq c_{\tau}$ for all $\bm{\xi} \in \mathbb{R}^{d}$ and a.e. $\bm{x}  \in \Gamma_C$. 

\medskip
\noindent
To obtain a weak formulation of Problem~$P$ we consider the following Hilbert spaces
\begin{align*}
 &\mathcal{H} = L^2(\Omega;\mathbb{S}^{d}), \qquad
 V = \{\bm{v} \in H^1(\Omega)^d\ |\ \bm{v} = \bm{0} \textrm{ on } \Gamma_{D}\},
\end{align*}
endowed with the inner scalar products
\begin{align*}
 &(\bm{\sigma},\bm{\tau})_{\mathcal{H}} = \int_\Omega \sigma_{ij}\tau_{ij} \, dx, \qquad
 (\bm{u},\bm{v})_V = (\bm{\varepsilon}(\bm{u}),\bm{\varepsilon}(\bm{v}))_{\mathcal{H}},
\end{align*}
respectively.
The fact that space $V$ equipped with the corresponding norm $\|\cdot\|_V$ is complete follows from Korn's inequality, and its application is allowed because we assume that $meas(\Gamma_{D}) > 0$.
We consider the trace operator $\gamma \colon V \to L^2(\Gamma_{C})^d=X$.

Using standard procedure, the Green formula and the definition of generalized subdifferential, we obtain a weak formulation of Problem~$P$ in the form of hemivariational inequality.

\medskip
\noindent
\textbf{Problem $\bm{P_{hvi}}$:} {\it \ Find a displacement $\bm{u} \in V$ such that for all $\bm{v} \in V$}
\begin{align}
  &\langle A\bm{u}, \bm{v} \rangle_{V^*\times V}  +\int_{\Gamma_C} j_2^0 (\bm{x}, \gamma \bm{u}(\bm{x})); \gamma \bm{v}(\bm{x}))\, da \geq \langle \bm{f}, \bm{v} \rangle_{V^*\times V}. \label{PV}
\end{align}

\medskip
\noindent
Here, the operator $A \colon V \to V^*$ and $\bm{f} \in V^*$ are defined for all $\bm{w},\bm{v} \in V$ as follows
\begin{align}
 &\langle A\bm{w}, \bm{v} \rangle_{V^*\times V} = (\mathcal{A}(\bm{\varepsilon}(\bm{w})),\bm{\varepsilon}(\bm{v}))_{\mathcal{H}}, \label{defA}\\
 &\langle \bm{f}, \bm{v} \rangle_{V^* \times V} = \int_{\Omega}\bm{f}_{0}\cdot \bm{v}\, dx + \int_{\Gamma_{N}}\bm{f}_{N}\cdot \gamma \bm{v}\, da \label{defF}
\end{align}
and $j\colon \Gamma_C\times \mathbb{R}^d \to \mathbb{R}$ is defined for all $\bm{\xi} \in \mathbb{R}^d$ and $\bm{x}\in \Gamma_C$ by
\begin{align}
j(\bm{x}, \bm{\xi}) =  j_{\nu}(\bm{x}, \xi_{\nu}) + h_{\tau}(\bm{x})\, j_{\tau}(\bm{x}, \bm{\xi}_{\tau}). \label{Jj}
\end{align}
It is easy to check that under assumptions $H(\mathcal{A})$ and $(H_0)$ and by the Sobolev trace theorem the operator~$A$, the functional~$\bm{f}$ and $\gamma \in  \mathcal{L}(V, X)$ satisfy $H(A)$ and \ref{assumption:f}, respectively.
We also define the functional $J \colon L^2(\Gamma_C)^d \to \mathbb{R}$ for all $\bm{v} \in L^2(\Gamma_C)^d $ by
\begin{align}
J(\bm{v}) =  \int_{\Gamma_C} j(\bm{x}, \bm{v}(\bm{x}))\, da,\label{J}
\end{align}

\medskip
\noindent
We remark that functional $J$ defined by {\rm{(\ref{Jj})-(\ref{J})}} under assumptions $H(j_{\nu})$, $H(j_{\tau})$ and $H(h)$ satisfies $H(J)$ (cf. \cite[Lemma 4]{JO}).

\medskip
\noindent
With the above properties, we have the following existence and uniqueness result for Problem~$P_{hvi}$.
\begin{theorem} \label{O=H}
If assumptions $H(\mathcal{A})$, $H(j_{\nu})$, $H(j_{\tau})$, $H(h)$, $(H_0)$ and $(H_s)$ hold, then Problems~$P_{hvi}$ and $P_{incl}$ (with functional $J$ dependent only on one variable) are equivalent. Moreover, they have a unique solution~$\bm{u}\in V$ and this solution satisfies
\begin{align*}
  &\| \bm{u} \|_V \leq c\, (1+\| \bm{f} \|_{V^*})
\end{align*}
with a positive constant $c$.
\end{theorem}

\begin{pfs}
We notice that the assumptions of Theorem~\ref{I=O} are satisfied. This implies that Problem~$P_{incl}$ has a unique solution. If $\bm{u}\in V$ is a solution to Problem $P_{incl}$ then it satisfies $\langle \bm{f} - A\bm{u}, \bm{v} \rangle_{V^*\times V} \leq J^0(\gamma \bm{u}; \gamma \bm{v})$ for all $\bm{v}\in V$. Hence, by Corollary~4.15~(iii) in \cite{MOS} we get that every solution to Problem $P_{incl}$ solves Problem $P_{hvi}$. Using a similar technique as in the proof of Theorem~\ref{I=O}, we can show that if Problem $P_{hvi}$  has a solution, it is unique. Combining these facts we obtain our assertion.
\end{pfs}

\section{Methods overview} \label{sec:mo}

In this section we present a brief overview of three established algorithms for solving contact problems - direct optimization method, augmented Lagrangian method and primal-dual active set strategy. References are also provided for more detailed treatment of each method. In versions presented here, all listed algorithms employ Finite Element Method (FEM). Let us start with the first mentioned method.

\subsection{Direct optimization method}

The idea behind direct optimization method is to replace weak formulation of contact problem with equivalent minimization problem. In the case of Problem~$P_{hvi}$ it takes the form

\medskip
\noindent
\textbf{Problem $\bm{\widehat{P}_{opt}^{h}}$:} {\it \ Find  $\bm{u}^h \in V^h$ such that}
\begin{align*}
0 \in \partial \widehat{\mathcal{L}}(\bm{u}^h),
\end{align*}
where functional $\widehat{\mathcal{L}}: V \rightarrow \mathbb{R}$ is defined for all $\bm{v} \in V$ as follows
\begin{align*}
  \widehat{\mathcal{L}}(\bm{v}) = \frac{1}{2} \langle A  \bm{v} ,  \bm{v}  \rangle_{V^*\times V}  - \langle  f ,  \bm{v}  \rangle_{V^*\times V} + J(\gamma  \bm{v}  ),
\end{align*}
and operators $A$, $\bm{f}$ and $J$ are defined by~(\ref{defA}),~(\ref{defF}) and~(\ref{J}), respectively.

Even though functions $j_{\nu}$ and $j_{\tau}$ can be nonmonotone and nonconvex, because of relaxed monotonicity condition on their subdifferentials combined with smallness assumption, $\widehat{\mathcal{L}}$ is still convex. Direct optimization method in a more complex setting (with function $h_{\tau}$ dependent on $u_\nu(\bm{x})$) using the Uzawa algorithm is presented in \cite{JO}.


\subsection{Augmented Lagrangian method}
Augmented Lagrangian method is a technique which regularizes nondifferentiable terms governing body behavior on contact boundary by addition of auxiliary Lagrange multipliers. These multipliers, represented by $\bm{\lambda^h}$, can be interpreted as normal and tangential forces acting on the body. Augmented Lagrangian approach expresses Problem~$P_{hvi}$ as a system of nonlinear equations.

Let us introduce $\Gamma_C^h$ as a discretization of contact interface based on FEM mesh, such that it consists of those nodes and edges of triangles that represent boundary $\Gamma_C$. Let $n_{C}$ be the number of independent points on boundary $\Gamma_C^h$ and $n_{tot}$ be the total number of independent nodes (outside of boundary $\Gamma_D$) on the FEM mesh.

In order to introduce the space of Lagrange multipliers, we use a contact element composed by one edge of $\Gamma_C^h$ and one Lagrange multiplier node. In our case FEM with affine polynomials is used for the displacement and FEM with constant polynomials is used for multipliers. This leads to space $H^h_{\Gamma_C}$, containing linear combinations of piecewise constant functions equal to $1$ on one edge of $\Gamma_C^h$ and $0$ everywhere else.

Let us now introduce matrices $W^h = (w_{ij})\in \mathbb{R}^{d \times n_{tot}}$ and $C^h = (c_{ij}) \in \mathbb{R}^{d \times n_{C}}$ with $w_{ij}$ and $c_{ij}$ being $i$-th coordinates of functions $\bm{v}^h \in V^h$ and $\bm{\gamma}^h \in H^h_{\Gamma_C}$ at node $j$ of FEM mesh, respectively (recall that $d$ is the dimension of considered body). We now reshape those matrices to obtain vectors $\widehat{\bm{v}}^h \in \mathbb{R}^{d \cdot n_{tot}}$ and $\widehat{\bm{\gamma}}^h \in \mathbb{R}^{d \cdot n_{C}}$, containing consecutive rows of respective matrices stacked in sequence. The generalized elastic term $\widetilde{G}(\widehat{\bm{u}}^h) \in \mathbb{R}^{d \cdot n_{tot}} \times \mathbb{R}^{d \cdot n_{C}}$ is defined by $\widetilde{G}(\widehat{\bm{u}}^h) = (G(\widehat{\bm{u}}^h), \bm{0}_{\mathbb{R}^{d \cdot n_{C}}})$, where $\bm{0}_{\mathbb{R}^{d \cdot n_{C}}}$ is the zero element of $\mathbb{R}^{d \cdot n_{C}}$ and $G(\widehat{\bm{u}}^h)$ denotes the term given for all $\bm{u}^h, \bm{v}^h \in V^h$ by
\begin{align*}
(G(\widehat{\bm{u}}^h)\cdot \widehat{\bm{v}}^h)_{\mathbb{R}^{d\cdot n_{tot}}} = \langle A  \bm{u}^h  -\bm{f},  \bm{v}^h  \rangle_{V^*\times V}.
\end{align*}
Here, operators $A$ and $\bm{f}$ are defined by~(\ref{defA}),~(\ref{defF}), respectively. 
One possibility of contact operator $\mathcal{F}(\widehat{\bm{u}}^h, \widehat{\bm{\xi}}^h)$, that deals with contact effects, can be defined for all $\widehat{\bm{u}}^h, \widehat{\bm{v}}^h \in \mathbb{R}^{d\cdot n_{tot}}, \widehat{\bm{\lambda}}^h, \widehat{\bm{\gamma}}^h \in \mathbb{R}^{d\cdot n_{C}}$, $\bm{u}^h, \bm{v}^h \in V^h$ and $\bm{\lambda}^h, \bm{\gamma}^h \in H^h_{\Gamma_C}$ by
\begin{align*}
&(\mathcal{F}(\widehat{\bm{u}}^h, \widehat{\bm{\lambda}}^h) \cdot (\widehat{\bm{v}}^h, \widehat{\bm{\gamma}}^h))_{\mathbb{R}^{d\cdot n_{tot} + d\cdot n_{C}}}\\
&\qquad + \int_{\Gamma_C} \, p(u^h_{\nu})\cdot \bm{v}^h da\\
&\qquad = \int_{\Gamma_C} \, \nabla_{\bm{w}} [l_{\nu}(\bm{u}^h, \bm{\lambda}^h) + l_{\tau}(\bm{u}^h, \bm{\lambda}^h)] \cdot \bm{v}^h da\\
&\qquad + \int_{\Gamma_C} \, \nabla_{\bm{\lambda}} [l_{\nu}(\bm{u}^h, \bm{\lambda}^h) + l_{\tau}(\bm{u}^h, \bm{\lambda}^h)] \cdot \bm{\gamma}^h da,
\end{align*}
where $\nabla_{\bm{y}}$ represents the gradient operator with respect to the variable $\bm{y}$. Functions $p$, $l_{\nu}$ and $l_{\tau}$ depend on the specific problem and are defined in Section~\ref{sec:sim}. The augmented Lagrangian approach is now expressed by the following system of equations.

\medskip
\noindent
\textbf{Problem $\bm{P_{Lag}^{h}}$:} {\it \ Find a displacement $\widehat{\bm{u}}^h \in \mathbb{R}^{d \cdot n_{tot}}$ and a stress multiplier field $\widehat{\bm{\lambda}}^h \in \mathbb{R}^{d \cdot n_{C}}$ such that}
\begin{align*}
\widetilde{G}(\widehat{\bm{u}}^h) + \mathcal{F}(\widehat{\bm{u}}^h, \widehat{\bm{\lambda}}^h) = \bm{0}.
\end{align*}
Note that using this method we approximate not only value of $\bm{u}$, but also values $\sigma_{\nu}$ and $\bm{\sigma}_{\tau}$ by $\widehat{\lambda}_{\nu}^h$ and $\widehat{\bm{\lambda}}_{\tau}^h$, respectively. More complete description of the presented approach can be found in \cite{BBK, BBKR}. For further details about discretization of contact interface and augmented Lagrangian method in general we refer to \cite{AC, PC, W}.

\subsection{Primal-dual active set strategy}
In primal-dual active set strategy we keep track of all points on discretized boundary $\Gamma^h_C$, and assign them to sets that reflect different ``parts'' of boundary laws $\partial j_\nu$ and $\partial j_\tau$. For example, every active set for $\partial j_\tau$ corresponds to a single or multivalued part of this function and can be interpreted as the physical state of the point assigned to this set (e.g. stick vs. slip zone). The division into these sets is specific to each chosen function and is reflected in implementation.

The main idea behind this strategy is to simplify integral over contact boundary in Problem $P_{hvi}$. We can do that using conditions that are implied by assignment of any point $\bm{x} \in \Gamma_C$ to a specific set (under the assumption that this assignment is correct). Initially all points are assigned to sets corresponding to $u_\nu =0$ and $\bm{u}_{\tau}=\bm{0}$. In the iterative procedure we successively solve the simplified problem with points assigned to selection of sets, and reassign them after each iteration. This  reassignment is conducted according to specific rules, so that we can move only between adjacent ``parts" of the graph of functions $\partial j_\nu$ and $\partial j_\tau$. We repeat this procedure until convergence to finally obtain a solution with all points in correct sets with respect both to $\partial j_\nu$ and $\partial j_\tau$. Further details of primal-dual active set strategy can be found in \cite{K}, and applications of this algorithm are presented for example in \cite{BCKYZ, XCHQ}.

\section{Simulations} \label{sec:sim}

We now consider Problem~$P$ and its weak formulation  Problem~$P_{hvi}$. We use data within previously presented theoretical framework, so the resulting model satisfies required assumptions and can be easily implemented using all three introduced methods.

\subsection{Data}

We set $d=2$ and consider a rectangular set $\Omega = [0,1] \times [0,1]$ with following partition of the boundary
\begin{align*}
 &\Gamma_{D} = \{0\} \times [0,1], \quad \Gamma_{N} = ([0,1] \times \{1\}) \cup (\{1\} \times [0,1]), \quad \Gamma_{C} = [0,1] \times \{0\}.
\end{align*}
The elasticity operator $\mathcal{A}$ is defined by
\begin{align*}
 &\mathcal{A}(\bm{\tau}) = 2\eta\bm{\tau} + \lambda \mbox{tr}(\bm{\tau})I,\qquad \bm{\tau} \in \mathbb{S}^2.
\end{align*}
Here, $I$ is the identity matrix, $\mbox{tr}$ denotes the trace of the matrix, $\lambda$ and $\eta$ are the Lam\'e coefficients, $\lambda, \eta >0$.
In our simulations we take the following data
\begin{align*}
 &\lambda = \eta = 4, \\
 &\bm{u}_{0}(\bm{x}) = (0,0), \quad \bm{x} \in \Omega,\\
 &\bm{f}_N(\bm{x}) = (0,0), \quad \bm{x} \in \Gamma_N.
\end{align*}
 We take nondifferentiable and nonconvex function $j_{\nu}$ and nondifferentiable function $j_{\tau}$ such that
\begin{align*}
 j_\nu(\bm{x}, \xi) &= \left \{ \begin{array}{ll}
   0, &\xi \in (-\infty,\, 0), \\
   \frac{p_{const}}{2}\, \xi^2 + q_{max}\, \xi , &\xi \in [0,\, \infty), \\
  \end{array} \right. \bm{x} \in \Gamma_C,\\[2mm]
 j_{\tau}(\bm{x}, \bm{\xi}) &= \|\bm{\xi}\|, \quad \bm{\xi} \in \mathbb{R}^2,\ \bm{x} \in \Gamma_C,
\end{align*}
where $q_{max}, p_{const} > 0$. This choice corresponds to (\ref{P5}) and (\ref{P6}) with
\begin{align*}
 \partial j_{\nu}(\bm{x}, \xi) &= p(\bm{x}, \xi) + \partial q(\bm{x}, \xi), \quad \xi \in \mathbb{R},\ \bm{x} \in \Gamma_C,\\[2mm]
 \partial j_{\tau}(\bm{x}, \bm{\xi}) &= \left \{ \begin{array}{ll}
   -1, &\xi \in (-\infty,\, 0), \\
   \mbox{[} {-1},\ 1\mbox{]}, &\xi = 0, \\
   1, &\xi \in (0,\, \infty),
  \end{array} \right. \bm{x} \in \Gamma_C,
 \end{align*}
where
\begin{align*}
 p(\bm{x}, \xi) &= \left \{ 
 \begin{array}{ll}
   0, &\xi \in (-\infty,\, 0), \\
   p_{const}\, \xi, &\xi \in [0,\, \infty), \\
 \end{array} \right. \bm{x} \in \Gamma_C,\\[2mm]
 \partial q(\bm{x}, \xi) &= \left \{ 
 \begin{array}{ll}
   0, &\xi \in (-\infty,\, 0), \\
       \mbox{[} 0,\ q_{max}\mbox{]}, &\xi = 0, \\
   q_{max}, &\xi \in (0,\, \infty), \\
 \end{array} \right. \bm{x} \in \Gamma_C.
\end{align*}

Remaining data is introduced for each example later, so that by applying various modifications of input parameters we can observe change in reaction of the body. It can be checked that all selected functions and data satisfy corresponding assumptions \ref{assumption:$A$}, \ref{assumption:j_nu}, \ref{assumption:j_tau}, \ref{assumption:h}, \ref{assumption:H_0}. Hence, we know that Problem~$P_{hvi}$ has a unique solution that can be estimated numerically.

\subsection{Implementation details}

We employ FEM and use space $V^h$ of continuous piecewise affine functions as a family of approximating subspaces. Uniform discretization of the problem domain according to the spatial discretization parameter $h$ is used. The contact boundary $\Gamma_C$ is divided into $1/h$ equal parts. For each example we start simulations with $h = 1/2$, which is successively halved. For all presented algorithms we can choose starting point $\bm{u}_0^h$ arbitrarily, but this choice affects convergence. To speed it up, we use ``warm start'' procedure: for mesh size $h=1/2$ we take $\bm{u}_0^h=\bm{0}$, and for every other $h$ we use a solution to problem with mesh size $h/2$.

To decrease the dimension of the considered discrete problem, we use the Schur complement method described in \cite{P}. This method reformulates a discrete scheme, so that we have to solve only for vertices located on the contact boundary $\Gamma_C$. Then we can retrieve the solution for all other vertices using simple matrix inversion and multiplication. These operations still  consume a lot of resources, but are feasible for much bigger mesh size $h$.

{
The code is written in Python, partly using Cython to obtain better performance. The solution is calculated using implementations of solvers from SciPy library. For augmented Lagrangian method we use {\it fsolve} to solve $\bm{\widehat{P}_{Lag}^{h}}$, and for primal-dual active set strategy we use {\it newton\_krylov} to solve discretized and simplified version of $P_{hvi}$. For direct optimization method we use {\it minimize} function with Powell's conjugate direction method to solve $\bm{\widehat{P}_{opt}^{h}}$, because minimized functional $\widehat{\mathcal{L}}$ is not necessarily differentiable. We can also employ other nonsmooth optimization algorithms such as the proximal bundle method  (see \cite{BKM}) or methods empirically proven to work for nondifferentiable functions, such as Broyden–Fletcher–Goldfarb–Shanno algorithm. Each presented implementation is chosen as one that gives best performance for each method. In all cases, we used default stopping criteria, based on change of argument value and change of function value with a maximal limit of iterations.
}

Let us now present how chosen data correspond to implementation details for augmented Lagrangian method and primal-dual active set strategy.

\medskip
For the augmented Lagrangian method, we introduce positive penalty coefficients $\epsilon_\nu, \epsilon_\tau > 0$. Definitions of functions $j_{\nu}$ and $j_{\tau}$ translate to auxiliary operators $l_{\nu}$ and $l_{\tau}$ in the form of
\begin{align*}
&l_{\nu}(\bm{u}^h, \bm{\lambda}^h) = \left \{ \begin{array}{ll}
   q_{max}\, u^h_\nu, &\lambda^h_\nu + \epsilon_\nu u^h_{\nu} \in (-\infty,\, -q_{max}], \\[2mm]
    \lambda^h_\nu\, u^h_\nu + \frac{\epsilon_\nu}{2}\, (u^h_\nu)^2, &\lambda^h_\nu + \epsilon_\nu u^h_{\nu} \in (-q_{max},\, 0], \\[2mm]
   -\frac{1}{2\epsilon_\nu}\, (\lambda_\nu)^2, &\lambda^h_\nu + \epsilon_\nu u^h_{\nu}\in (0,\, \infty),
  \end{array} \right.\\[3mm]
&l_{\tau}(\bm{u}^h, \bm{\lambda}^h) = \left \{ \begin{array}{ll}
   \bm{\lambda}^h_{\bm{\tau}}\cdot \bm{u}^h_\tau + \frac{\epsilon_\tau}{2}\, \bm{u}^h_\tau \cdot \bm{u}^h_\tau, &\|\bm{\lambda}^h_{\bm{\tau}} + \epsilon_\tau \bm{u}^h_{\tau}\| \leq h_\tau, \\[2mm]
   -\frac{1}{2\epsilon_\tau}\, \big[ \|\bm{\lambda}^h_{\bm{\tau}}\|^2 - 2 h_\tau \|\bm{\lambda}^h_{\bm{\tau}} \| + (h_\tau)^2 \big], &\|\bm{\lambda}^h_{\bm{\tau}} + \epsilon_\tau \bm{u}^h_{\tau}\| > h_\tau.
  \end{array} \right.
\end{align*}
In the iterative procedure we solve system of equations present in $P_{Lag}^{h}$ for fixed values of penalty parameters $\epsilon_\nu, \epsilon_\tau$, decrease them, and, taking previously obtained solution as a starting point, repeat until convergence. 

\bigskip
In the case of primal-dual active set strategy, function $\partial j_\nu$  divides points of $\Gamma_C$ into exclusive sets $N^1$, $N^2$, $N^3$ and function $\partial j_\tau$ into exclusive sets $T^1$, $T^2$. For any point $\bm{x} \in \Gamma^h_C$ we have the following possibilities
\begin{itemize}
  \item if $\bm{x} \in N^1$ then $u^h_{\nu}(\bm{x}) < 0$, which implies $\sigma^h_\nu(\bm{x}) = 0$ (points on the boundary lifted from the foundation),\itemsep = 0.3em
  \item if $\bm{x} \in N^2$ then $u^h_{\nu}(\bm{x}) = 0$, which implies $-\sigma^h_\nu(\bm{x}) \in \partial q(\bm{x}, 0)$ (points in contact with the foundation experiencing force in normal direction below or equal to specified threshold $q_{max}$, i.e. in rigid state),\itemsep = 0.3em
  \item if $\bm{x} \in N^3$ then $u^h_{\nu}(\bm{x}) > 0$, which implies $-\sigma^h_\nu(\bm{x}) = p_{const}\, \xi + q_{max}$ (points as described above, but with force over specified threshold $q_{max}$, i.e. in flexible state),\itemsep = 0.7em
  \item if $\bm{x} \in T^1$ then $\|\bm{u}^h_{\tau}(\bm{x})\| = 0$, which implies $-\bm{\sigma^h}_\tau(\bm{x}) \in [-h_\tau,\, h_\tau] $ (points on the boundary experiencing force in tangential direction with norm below or equal to friction bound $h_\tau$, i.e. in the stick zone),\itemsep = 0.3em
  \item if $\bm{x} \in T^2$ then $\|\bm{u}^h_{\tau}(\bm{x})\| > 0$, which implies $-\bm{\sigma^h}_\tau(\bm{x}) \in \{-h_\tau,\, h_\tau\} $ (points as described above, but with norm of force over friction bound $h_\tau$, i.e. in the slip zone).\itemsep = 0.3em
\end{itemize}
Initially all points are assigned to $N^1_0$ and $T^1_0$, with subscript denoting current iteration. Then the following rules (with $\epsilon > 0$ being a small value added for numerical stability) are applied
\begin{itemize}
  \item if $\bm{x}\in N^1_i$ and $u^h_\nu(\bm{x}) \in (-\infty, -\epsilon)$, then $\bm{x}\in N^1_{i+1}$,\itemsep = 0.3em
  \item if $\bm{x}\in N^1_i$ and $u^h_\nu(\bm{x}) \in [-\epsilon, \infty)$, then $\bm{x}\in N^2_{i+1}$,\itemsep = 0.7em
  \item if $\bm{x}\in N^2_i$ and $\sigma^h_\nu(\bm{x}) \in (-\infty, - \epsilon)$, then $\bm{x}\in N^1_{i+1}$,\itemsep = 0.3em
  \item if $\bm{x}\in N^2_i$ and $\sigma^h_\nu(\bm{x}) \in [-\epsilon, q_{max} + \epsilon)$, then $\bm{x}\in N^2_{i+1}$,\itemsep = 0.3em
  \item if $\bm{x}\in N^2_i$ and $\sigma^h_\nu(\bm{x}) \in [q_{max} + \epsilon, \infty)$, then $\bm{x}\in N^3_{i+1}$,\itemsep = 0.7em
  \item if $\bm{x}\in N^3_i$ and $u^h_\nu(\bm{x}) \in (-\infty, -\epsilon)$, then $\bm{x}\in N^2_{i+1}$,\itemsep = 0.3em
  \item if $\bm{x}\in N^3_i$ and $u^h_\nu(\bm{x}) \in [-\epsilon, \infty)$, then $\bm{x}\in N^3_{i+1}$,\itemsep = 0.7em
  \item if $\bm{x}\in T^1_i \cup T^2_i$ and $\|\bm{\sigma}^h_\tau(\bm{x})\| < h_\tau + \epsilon$, then $\bm{x}\in T^1_{i+1}$,\itemsep = 0.3em
  \item if $\bm{x}\in  T^1_i \cup T^2_i$, and $\|\bm{\sigma}^h_\tau(\bm{x})\| \geq h_\tau + \epsilon$, then $\bm{x}\in T^2_{i+1}$.
\end{itemize}

Estimation of $\bm{\sigma}^h(\bm{x})$ required for this step can be calculated from discretization of constitutive law, using values of displacement $\bm{u}^h$ of neighbors of $\bm{x}$ on the FEM mesh.

\subsection{ Final results}

Finally, we present outputs obtained in our simulations and report empirical estimation of numerical errors. In examples we choose mesh corresponding to $h=1/32$ for better visibility. As expected, all considered algorithms give similar solutions, so for each example we select the output of one algorithm for illustration. We plot deformation of the body and forces acting on the contact interface (mirrored with respect to the boundary for visibility).

\begin{figure}[h!]
\begin{minipage}{.5\textwidth}
Figure~\ref{fig:lagrangian} presents result for data
\begin{equation*}
 \begin{split}
 &\bm{f}_0(\bm{x}) = (-0.5,\, -1.0), \quad \bm{x} \in \Omega,\\
 &h_{\tau} = 0.1, \quad q_{max} = 0.1, \quad p_{const} = 10.
 \end{split}
\end{equation*}

\medskip
We push the body down and to the left with force $\bm{f}_0$. In this case the coefficient $p_{const}$ has the highest influence on response of the foundation in normal direction. It causes forces to increase gradually with penetration and models a foundation made of a soft material. A small influence of friction can also be observed.
\end{minipage}
\begin{minipage}{.5\textwidth}
  \centering
    \includegraphics[width=0.9\linewidth]{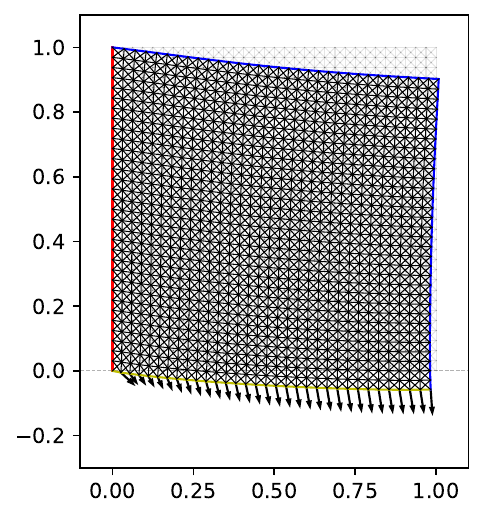}
    \caption{Output obtained\\ using augmented Lagrangian method}
    \label{fig:lagrangian}
\vspace{10mm}
\end{minipage}
\end{figure}

\begin{figure}[h!]
\begin{minipage}{.5\textwidth}
Figure~\ref{fig:opt} presents result for data
\begin{equation*}
 \begin{split}
 &\bm{f}_0(\bm{x}) = (-0.5,\, -1.0), \quad \bm{x} \in \Omega,\\
 &h_{\tau} = 0.1, \quad q_{max} = 0.7, \quad p_{const} = 0.
 \end{split}
\end{equation*}

\medskip
Here forces in normal direction increase up to a factor $q_{max}$ and, because $p_{const}=0$, stop increasing any further. The foundation response is therefore limited by this factor.
\end{minipage}
\begin{minipage}{.5\textwidth}
  \centering
    \includegraphics[width=0.9\linewidth]{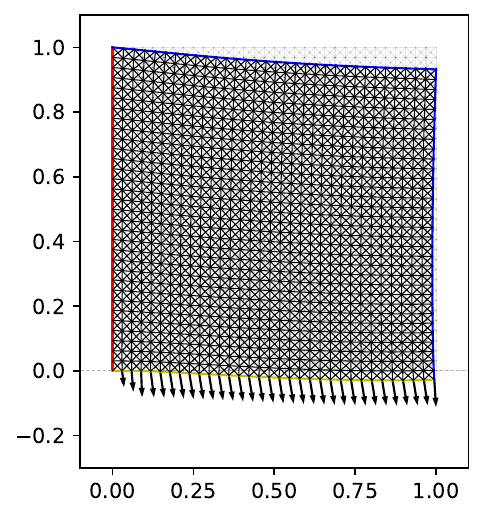}
    \caption{Output obtained\\using direct optimization method} \label{fig:opt}
\vspace{10mm}
\end{minipage}
\end{figure}

\begin{figure}[h!]
\begin{minipage}{.5\textwidth}
Figure~\ref{fig:primal-dual} presents result for data
\begin{equation*}
 \begin{split}
 &\bm{f}_0(\bm{x}) = (-0.5,\, -1.0), \quad \bm{x} \in \Omega,\\
 &h_{\tau} = 0.5, \quad q_{max} = 0.5, \quad p_{const} = 0.
 \end{split}
\end{equation*}

\bigskip
In this example we increase the friction bound $h_\tau$. As we push the body to the left, points on the left side of the boundary move to the slip zone while points on the right side cannot overcome friction bound and stay in the stick zone.
\end{minipage}
\begin{minipage}{.5\textwidth}
  \centering
    \includegraphics[width=0.9\linewidth]{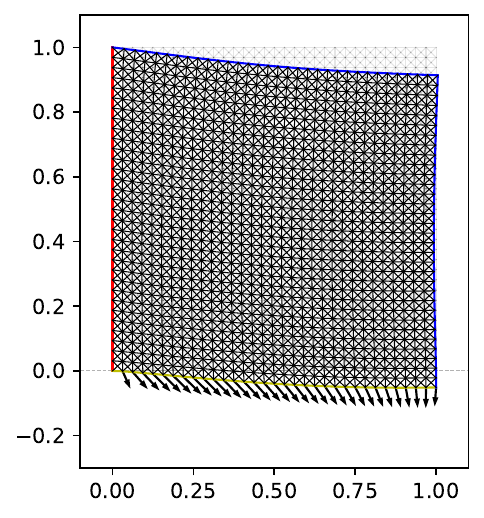}
    \caption{Output obtained\\using primal-dual active set strategy} \label{fig:primal-dual}
\vspace{10mm}
\end{minipage}
\end{figure}

\begin{figure}[h!]
\begin{minipage}{.5\textwidth}
Figure~\ref{fig:opt2} presents result for data
\begin{equation*}
 \begin{split}
 &\bm{f}_0(\bm{x}) = (0.5,\, -1.0), \quad \bm{x} \in \Omega,\\
 &h_{\tau} = 0.1, \quad q_{max} = 10, \quad p_{const} = 0.
 \end{split}
\end{equation*}

\medskip
We change force $\bm{f}_0$ and, as a result, the body is pushed to the right. We also set $q_{max} = 10$, and this effectively enforces Signorini condition. We can see that the body cannot penetrate the foundation and we can also observe the influence of friction forces.
\end{minipage}
\begin{minipage}{.5\textwidth}
  \centering
    \includegraphics[width=0.9\linewidth]{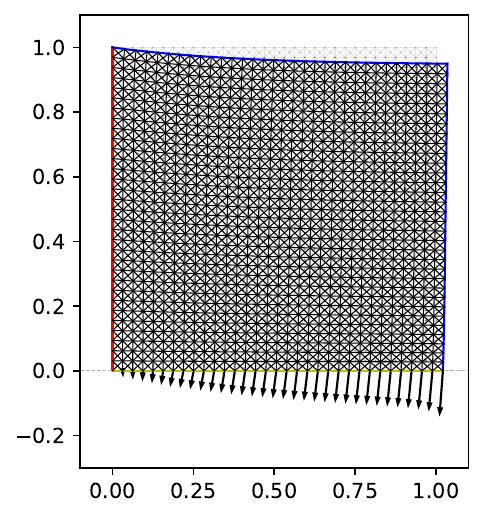}
    \caption{Output obtained\\using direct optimization method} \label{fig:opt2}
\vspace{10mm}
\end{minipage}
\end{figure}

\begin{figure}[h!]
\centering
  \includegraphics[width=0.6\linewidth]{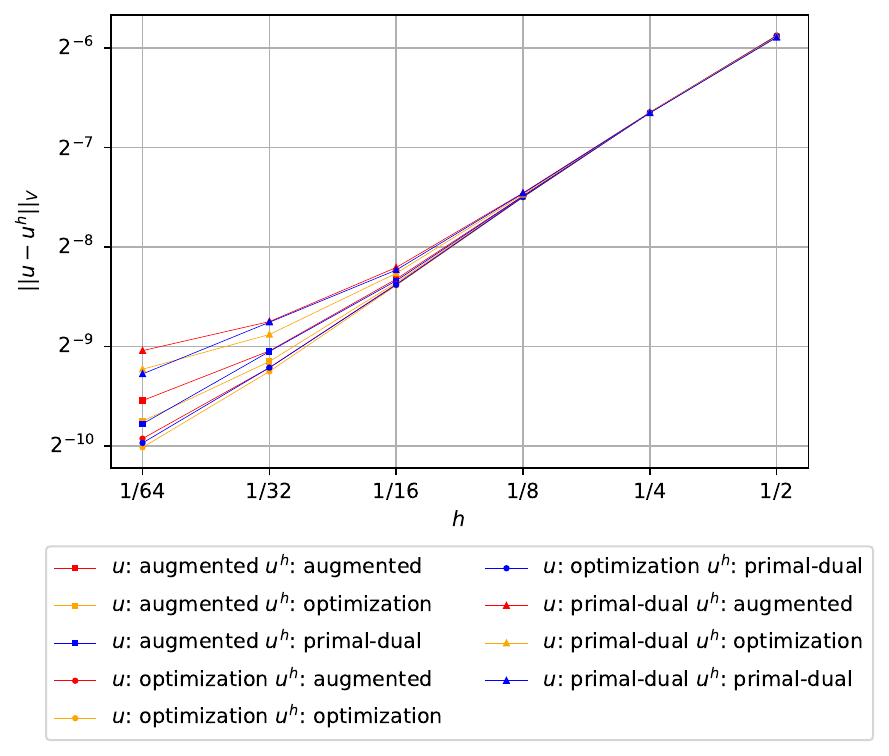}
  \caption{Numerical errors} \label{fig:errors}
\end{figure}

A comparison of numerical errors $\|\bm{u} - \bm{u}^h\|_V$ computed for a sequence of solutions to discretized problems on a model problem with data
\begin{equation}\label{e5}
 \begin{split}
 &\bm{f}_0(\bm{x}) = (-0.8,\, -0.8), \quad \bm{x} \in \Omega,\\
 &h_{\tau} = 0.5, \quad q_{max} = 0.3, \quad p_{const} = 2
 \end{split}
\end{equation}
is presented in Figure \ref{fig:errors}, where the dependence of the error estimate $\|\bm{u}  - \bm{u}^h\|_V$ with respect to $h$ is plotted on a log-log scale.
Because no analytical solution can be obtained, we took three numerical estimations with $h = 1/256$ and corresponding to each presented method as such ``exact'' solutions. All sequences of numerical solutions with varying $h$ were cross examined against each of ``exact'' solutions, giving 9 plots. We denote by $u$ ``exact" solutions (depending on chosen method), by $u^h$ sequence of numerical approximations (also for each method) and use abbreviations of presented methods' names. 
 {
 As we can see, in this case the primal-dual active set strategy and direct optimization method gave similar final estimations, closer to reference solutions than augmented Lagrangian method.

\begin{table}[ht]
\footnotesize
\centering
\begin{tabular}{ l r r r r r r }\hline
$h$ & & $1/8$ & $1/16$ & $1/32$ & $1/64$  & $1/128$  \\ \hline
Direct optimization & Time & $0.14s$ & $0.36s$ & $1.33s$ & $3.97s$ & $40.40s$  \\
& {\rm Functional evaluations} & 1027 & 1397 & 2832 & 2748 & 5329 \\ \hline
Augmented Lagrangian & Time & $0.26s$ & $0.74s$ & $1.38s$ & $3.79s$  &  $37.60s$ \\ \hline
Primal-Dual & Time & $0.07s$ & $0.19s$ & $0.76s$ & $6.96s$ & $178.07s$ \\
& {\rm Set iterations} & 4 & 4 & 4 & 5 &  6 \\ \hline
\end{tabular}
\caption{Computation time and number of iterations for each algorithm} \label{tabTwo}
\end{table}

In Table \ref{tabTwo} we summarised computation time and number of iterations for each algorithm. Presented results do not include time for computation of stiffness matrix, which is calculated beforehand and is the same for all methods. Additional metric "function evaluations" for direct optimization denotes how many times functional $\widehat{\mathcal{L}}$ was evaluated. "Set changes" for primal-dual denote how many iterations of assignments to sets $N$ and $T$ were performed. The fastest method for finer meshes in this case was augmented Lagrangian, closely followed by direct optimization method. Our implementation of primal-dual active set strategy was fastest for coarse, but slowest for fine mesh sizes.

We remark that direct optimization method was easiest to implement, followed by primal-dual and augmented Lagrangian, as it has most complicated interpretation. As stated before, augmented Lagrangian method simultaneously with $u$ on $\Omega$ calculates values of $\sigma$ on $\Gamma_C$, which for other methods had to be estimated from the value of $u$. We also remark that all presented results may vary depending on details of specific implementations.
 }

\bigskip
\noindent {\bf Acknowledgments}\\
The project leading to this application has received funding from the European Union's Horizon 2020 Research and Innovation Programme under the Marie Sklo\-do\-wska-Curie grant agreement no. 823731 CONMECH, from the Ministry of Science and Higher Education of Republic of Poland under Grant No 440328/PnH2/2019, and in part from National Science Centre, Poland under project OPUS no. 2021/41/B/ST1/01636.

\medskip
\noindent {\bf Conflict of interest}\\
The authors declare that they have no conflict of interest.



\end{document}